\documentclass[12pt]{article}
\usepackage[final]{epsfig}
\usepackage{graphics}
\usepackage{amsmath}
\usepackage{amsfonts}
\usepackage{latexsym}
\usepackage{amssymb}
\usepackage{graphicx}
\usepackage{epstopdf}
\newtheorem{lemma}{Lemma}

\newtheorem{theorem}[lemma]{Theorem}

\begin{document}
\newcommand{\eps}{{\varepsilon}}
\newcommand{\proofend}{$\Box$\bigskip}
\newcommand{\C}{{\mathbf C}}
\newcommand{\Q}{{\mathbf Q}}
\newcommand{\R}{{\mathbf R}}
\newcommand{\Z}{{\mathbf Z}}
\newcommand{\RP}{{\mathbf {RP}}}
\newcommand{\CP}{{\mathbf {CP}}}

\title {On algebraically integrable outer billiards}
\author{Serge Tabachnikov\\
{\it Department of Mathematics, Penn State University}\\
{\it University Park, PA 16802, USA}\\
e-mail: {\it tabachni@math.psu.edu}}
\date{\today}
\maketitle
\begin{abstract}
We prove that if the outer billiard map around a plane oval is algebraically integrable in a certain non-degenerate sense then the oval is an ellipse.
\end{abstract}

In this note, an outer billiard table is a compact convex domain in the plane bounded by an oval (closed smooth strictly convex curve) $C$. Pick a point $x$ outside of $C$. There are two tangent lines from $x$ to $C$; choose  one of them, say, the right one from the view-point of $x$, and reflect $x$ in the tangency point. One obtains a new point, $y$, and the transformation $T: x\mapsto y$ is the outer (a.k.a. dual) billiard map. We refer to \cite{D-T,Ta1,Ta2} for surveys of outer billiards.

If $C$ is an ellipse then the map $T$ possesses a 1-parameter family of invariant curves, the homothetic ellipses; these invariant curves foliate the exterior of $C$. Conjecturally, if an outer neighborhood of an oval $C$ is foliated by the invariant curves of the outer billiard map then $C$ is an ellipse -- this is an outer version of the famous Birkhoff conjecture concerning the conventional, inner billiards.  

In this note we show that ellipses are rigid in a much more restrictive sense of algebraically integrable outer billiards; see \cite{Bo} for the case of inner billiards.

We make the following assumptions. Let $f(x,y)$ be a (non-homogeneous) real polynomial such that zero is its non-singular value and $C$ is a component of the zero level curve. Thus $f$ is the defining polynomial of the curve $C$, and if a polynomial vanishes on $C$ then it is a multiple of $f$.
Assume that a neighborhood of $C$ is foliated by invariant curves of the outer billiard map $T$, and this foliation is algebraic in the sense that its leaves are components of the level curves of a real polynomial $F(x,y)$. Since $C$ itself is an invariant curve, we assume that $F(x,y)=0$ on $C$, and that $d F$ is not identically zero on $C$. Thus  $F(x,y)=g(x,y) f(x,y)$  where $g(x,y)$ is a polynomial, not identically zero on $C$. Under these assumptions, our result is as follows.

\begin{theorem} \label{main}
$C$ is an ellipse.
\end{theorem}

\paragraph{Proof.} Consider  the tangent vector field $v = F_y\ \partial/\partial x -F_x  \ \partial/\partial y$ (the symplectic gradient) along $C$. This vector field is non-zero (except, possibly, a finite number of points) and tangent to $C$. The tangent line to $C$ at point $(x,y)$ is given by $(x+\eps F_y,y-\eps F_x)$, and the condition that $F$ is $T$-invariant means that the function 
\begin{equation} \label{funct}
F(x+\eps F_y,y-\eps F_x)
\end{equation}
 is even in $\eps$ for  all $(x,y) \in C$. Expand in a series in $\eps$; the first non-trivial condition is cubic in $\eps$:
\begin{equation} \label{cubic}
W(F):=F_{xxx} F_y^3 -3 F_{xxy} F_y^2 F_x + 3 F_{xyy} F_y F_x^2 - F_{yyy} F_x^3 =0 \end{equation}
on $C$. We claim that this already implies that $C$ is an ellipse. The idea is that otherwise the complex curve $f=0$ would have an inflection point, in contradiction with identity (\ref{cubic}).

Consider the polynomial 
$$
H(F)= \det \left(
\begin{array}{cc}
F_y,& -F_x\\
F_{yy} F_x - F_{xy} F_y,& F_{xx} F_y- F_{xy} F_x 
\end{array}
\right).
$$

\begin{lemma} \label{lm}
One has:\\
1). $v(H(F))=W(F)$;\\
2). $H(F)=H(gf)=g^3 H(f)$ on $C$;\\
3). If $C'$ is a non-singular algebraic curve with a defining polynomial $f(x,y)$ and $(x,y)$ is an inflection point of $C'$ then $H(f)(x,y)=0$.
\end{lemma}

\paragraph{Proof of Lemma \ref{lm}.} The first two claims follow from straightforward computations. To prove the third, note that $H(f)$ is the second order term in $\eps$ of the Taylor expansion of the function $f(x+\eps f_y,y-\eps f_x)$ (cf.
 (\ref{funct})), hence $H(f)=0$ at an inflection point.
\proofend

It follows from Lemma \ref{lm}  and (\ref{cubic}) that $H(F)=$ const on  $C$. Since $C$ is convex, $H(F)\neq 0$, and we may assume that $H(F)=1$ on $C$. It follows that 
$g^3 H(f)-1$ vanishes on $C$ and hence
\begin{equation} \label{ideal}
g^3 H(f) - 1 = h f
\end{equation}
where $h(x,y)$ is some polynomial. 

Now consider the situation in $\CP^2$. We continue to use the notation $C$ for the complex algebraic curve given by the homogenized polynomial ${\bar f}(x:y:z)=f(x/z,y/z)$. Unless $C$ is a conic, this curve has inflection points (not necessarily real). Let $d$ be the degree of $C$.

\begin{lemma} \label{line}
Not all the inflections of $C$ lie on the line at infinity
\end{lemma}

\paragraph{Proof of Lemma \ref{line}.} Consider the Hessian curve given by
$$
\det \left(
\begin{array}{ccc}
{\bar f}_{xx},&{\bar f}_{xy},&{\bar f}_{xz}\\
{\bar f}_{yx},&{\bar f}_{yy},&{\bar f}_{yz}\\
{\bar f}_{zx},&{\bar f}_{zy},&{\bar f}_{zz}
\end{array}
\right)=0.
$$
The intersection points of the curve $C$ with its
Hessian curve are the inflection points of $C$ (recall that $C$ is non-singular). The degree of the Hessian curve is $3(d-2)$,  and by the Bezout theorem, the total number of inflections, counted with multiplicities, is $3d(d-2)$. Furthermore, the order of intersection equals the order of the respective inflection and does not exceed $d-2$, see, e.g, \cite{Wa}. The number of intersection points of $C$ with a line equals $d$, hence the inflection points of $C$ that lie on a fixed line contribute, at most, $d(d-2)$ to the total of $3d(d-2)$. The remaining inflection point lie off this line.
\proofend

To conclude the proof of Theorem \ref{main}, consider a finite inflection point of $C$.  
According to Lemma \ref{lm}, at such a point point, we have $f=H(f)=0$ which contradicts (\ref{ideal}). This is proves that $C$ is a conic.
\proofend

\paragraph {\bf Remarks.}
1). It would be interesting to remove the non-degeneracy assumptions in Theorem \ref{main}. 

2). A more general version of Birkhoff's integrability conjecture is as follows. Let $C$ be a plane oval whose outer neighborhood is foliated by closed curves. For a tangent line $\ell$ to $C$, the intersections with the leaves of the foliation define a local involution $\sigma$ on $\ell$. Assume that, for every tangent line, the  involution $\sigma$ is projective. Conjecturally, then $C$ is an ellipse and the foliation consists of ellipses that form a pencil (that is, share four -- real or complex -- common points). For a pencil of conics, the respective involutions are projective: this is a Desargues theorem, see \cite{Ber}. It
would be interesting to establish an algebraic version of this conjecture.
\bigskip

{\bf Acknowledgments}. Many thanks to Dan Genin for numerous stimulating conversations, to S. Bolotin for comments on his work \cite{Bo}, to V. Kharlamov for providing a proof of Lemma \ref{line} and to Rich Schwartz for interest and criticism.  The  author was partially supported by an NSF grant DMS-0555803.

\end{document}